\documentclass{article}
\usepackage{spconf,graphicx}

 \usepackage[english,french]{babel}   
  \usepackage{times}			
 
\usepackage[T1]{fontenc}
\usepackage[utf8]{inputenc} 
\DeclareUnicodeCharacter{00A0}{~}

\usepackage{amssymb,amsmath,amscd,amsfonts,amsthm,bbm,mathrsfs,yhmath}
\usepackage[shortlabels]{enumitem}

\DeclareMathOperator{\dom}{dom}

\DeclareMathOperator{\prox}{prox}

\DeclareMathOperator{\cl}{cl}

\newcommand{\bs}{\boldsymbol}

\newcommand{\sx}{{\mathsf x}}

\newcommand{\bP}{{{\mathbb P}}} 
\newcommand{\bE}{{{\mathbb E}}} 
\newcommand{\bN}{{{\mathbb N}}} 

\newcommand{\mC}{{\mathcal C}} 
\newcommand{\mD}{{\mathcal D}}

\newcommand{\sX}{{\mathsf X}}

\newcommand{\mcB}{{\mathscr B}}

\newcommand{\cG}{{{\mathcal G}}} 
 
\newcommand{\cD}{{{\mathcal D}}}

\newcommand{\cK}{{{\mathcal K}}}

\newcommand{\bR}{{{\mathbb R}}}

\newcommand{\ps}[1]{\langle #1 \rangle}

\usepackage[textwidth=2cm, textsize=footnotesize]{todonotes}  
\theoremstyle{definition}
\newtheorem{theorem}{Theorem}

\newtheorem{assumption}{Assumption}

\title{A Constant Step Stochastic Douglas-Rachford Algorithm with Application to Non Separable Regularizations}
\name{Adil Salim$^\star$, Pascal Bianchi$^\star$ and Walid Hachem$^\dagger$ \thanks{This work was supported by the Agence Nationale pour la Recherche,
France, (ODISSEE project, ANR-13-ASTR-0030) and by the Labex Digiteo-DigiCosme
(OPALE project), Universit\'e Paris-Saclay.
}}
\address{$\star$ LTCI, Télécom ParisTech, Université Paris-Saclay. \\
46, rue Barrault, 75634 Paris Cedex 13, France. \\
$\dagger$ CNRS / LIGM (UMR 8049), Universit\'e Paris-Est Marne-la-Vallée. 
\\
 5, boulevard Descartes, Champs-sur-Marne, 77454, Marne-la-Vallée Cedex 2, 
France.  }

\begin{document}
%
\maketitle
\begin{abstract}
The Douglas Rachford algorithm is an algorithm that converges to a minimizer of a sum of two convex functions. The algorithm consists in fixed point iterations involving computations of the proximity operators of the two functions separately. The paper investigates a stochastic version of the algorithm where both functions are random and the step size is constant. We establish that the iterates of the algorithm stay close to the set of solution with high probability when the step size is small enough. Application to structured regularization is considered.
\end{abstract}
\section{Introduction}
\label{sec:intro}
Many applications in the fields of machine learning~\cite{boyd2011distributed} and signal processing~\cite{chierchiaapproche} require the solution of the programming problem
\begin{equation}
\label{eq:minF+G}
\min_{x \in \sX} F(x) + G(x) 
\end{equation}
where $\sX$ is an Euclidean space, $F$ and $G$ are elements of the set $\Gamma_0(\sX)$ of convex, lower semi-continuous and proper functions. In these contexts, $F$ often represents a cost function and $G$ a regularization term. The Douglas-Rachford algorithm is one of the most popular approach towards solving Problem~\eqref{eq:minF+G}. Given $\gamma >0$, the algorithm is written 
\begin{align} 
\label{eq:dr-deter} 
y_{n+1} &= \prox_{\gamma F}(x_n) \nonumber\\
z_{n+1} & = \prox_{\gamma G}(2 y_{n+1} - x_n) \nonumber \\
x_{n+1} &= x_n + z_{n+1} - y_{n+1}
\end{align} 
where $\prox_{\gamma F}$ denotes the proximity operator of $F$, defined for every $x \in \sX$ by the equation 
\begin{equation*}
\prox_{\gamma F}(x) = \arg\min_{y \in \sX} \frac12 \|x - y\|^2 + \gamma F(y).  
\end{equation*}
Assuming that a standard qualification condition holds and that the set of solutions $\arg\min F+G$ of~\eqref{eq:minF+G} is not empty, the sequence $(y_n)_n$ converges to an element in $\arg\min F+G$ as $n \to +\infty$ (\cite{lio-mer-79,eckstein1992douglas}).

In this paper, we study the case where $F$ and $G$ are integral functionals of the form $$F({x}) = \bE_\xi(f({x},\xi)), \qquad G({x}) = \bE_\xi(g({x},\xi))$$ where $\xi$ is a random variable (r.v) from some probability space $(\Omega, \mathcal F, \bP)$ into a measurable space $(\Xi,\mathcal G)$, with distribution $\mu$, and where $\{f(\cdot,s), s \in \Xi\}$ and $\{g(\cdot,s), s \in \Xi\}$ are subsets of $\Gamma_0(\sX)$. In this context, the stochastic Douglas Rachford algorithm aims to solve Problem~\eqref{eq:minF+G} by iterating
\begin{align}
\label{eq:sto-dr}
  {y}_{n+1} &=\prox_{\gamma f(\cdot,\xi_{n+1})}({x}_{n}) \nonumber\\
  {z}_{n+1} &=\prox_{\gamma g(\cdot,\xi_{n+1})}(2{y}_{n+1}- {x}_{n})\nonumber\\
  {x}_{n+1} &= {x}_{n} + {z}_{n+1}-{y}_{n+1}\,,
\end{align}
where $(\xi_{n})_n$ is a sequence of i.i.d copies of the random variable $\xi$ and $\gamma > 0$ is the constant step size. Compared to the "deterministic" Douglas Rachford algorithm~\eqref{eq:dr-deter}, the stochastic Douglas Rachford algorithm~\eqref{eq:sto-dr} is an online method. The constant step size used make it implementable in adaptive signal processing or online machine learning contexts. In this algorithm, the function $F$ (resp. $G$) is replaced at each iteration $n$ by a random realization $f(\cdot,\xi_{n})$ (resp. $g(\cdot,\xi_{n})$). It can be implemented in the case where $F$ (resp. $G$) cannot be computed in its closed form~\cite{bia-hac-jota16,bia-hac-sal-(sub)jca17} or in the case where the computation of its proximity operator is demanding~\cite{sal-bia-hac-(sub)tac17}. Compared to other online optimization algorithm like the stochastic subgradient algorithm, the algorithm~\eqref{eq:sto-dr} benefits from the numerical stability of stochastic proximal methods.

Stochastic version of the Douglas Rachford algorithm have been considered in~\cite{chierchiaapproche,shi2013online}. These papers consider the case where $G$ is deterministic, \textit{i.e} is not written as an expectation and $F$ is written as an expectation that reduces to a sum. The latter case is also contained as a particular case of the algorithm~\cite{chambolle2017stochastic}. The algorithms~\cite{rosasco2015stochastic,combettes2016stochastic} are generalizations of a partially stochastic Douglas Rachford algorithm where $G$ is deterministic. The convergence of these algorithms is obtained under a summability assumption of the noise over the iterations. The stochastic Douglas Rachford studied in this paper was implemented in an adaptive signal processing context~\cite{mourya2017adaptive} to solve a target tracking problem.

Whereas the paper~\cite{mourya2017adaptive} is mainly focused on the application to target tracking, in this work we provide theoretical basis for the algorithm~\eqref{eq:sto-dr} and convergence results. Moreover, a novel application to solve a programming problem regularized with the overlapping group lasso online is provided.

The next section introduces some notations. Section~\ref{sec-th} is devoted to the statement of the main convergence result. In Section~\ref{sec-proof}, an outline of the proof of the result in Section~\ref{sec-th} is provided. Finally, the algorithm~\eqref{eq:sto-dr} is implemented to solve a regularized problem in Section~\ref{sec-simu}.

\section{Notations}

For every function $g \in \Gamma_0(\sX)$, $\partial g(x)$ denotes the subdifferential of $g$ at the point $x \in \sX$ and $\partial g_0(x)$ the least norm element in $\partial g(x)$. The domain of $g$ is denoted as $\dom(g)$. It is a known fact that the closure of $\dom(g)$, denoted as $\cl(\dom(g))$, is convex. For every closed convex set $\mC$, we denote by $\Pi_C$ the projection operator onto $\mC$. The indicator function of the set $\mC$ is defined by $\iota_{\mC}(x) = 0$ if $x \in \mC$, and $\iota_{\mC}(x) = +\infty$ elsewhere. It is easy to see that $\iota_{\mC} \in \Gamma_0(\sX)$ and that $\prox_{\iota_{\mC}} = \Pi_\mC$.

The Moreau envelope of $g \in \Gamma_0(\sX)$ is equal to
$$
g_\gamma({x}) = \min_{{y}\in \sX} g({y})+\frac{\|{y}-{x}\|^2}{2\gamma}
$$
for every ${x}\in \sX$. Recall that $g_\gamma$ is differentiable and $\nabla g_\gamma (x) = \frac{1}{\gamma}(x - \prox_{\gamma g}({x})).$ If $f \in \Gamma_0(\sX)$ is differentiable, then, $\partial f(x) = \{\nabla f(x)\}$ and $\nabla f(\prox_{\gamma f}(x)) = \nabla f_\gamma(x)$, for every $x \in \sX$.

When $S \subset \sX$, $d(x,S)$ denote the distance from the point $x \in \sX$ to the set $S$. In the context of algorithm~\eqref{eq:sto-dr} we shall denote $D(s) = \dom(g(\cdot,s))$ and $\cD = \dom(G)$. Denote $\mcB(\sX)$ the Borel sigma field over $\sX$. For every $p \geq 1$, $L^p(\Xi,\sX)$ is the set of all r.v $\varphi$ from the probability space $(\Xi,\cG,\mu)$ into the measurable space $(\sX,\mcB(\sX))$, such that $\|\varphi\|^p$ is integrable.
 
From now on, we shall state explicitly the dependence of the iterates of the algorithm in the step size and the starting point. Namely, we shall denote $(x_n^{\gamma,\nu})_n$ the sequence $(x_n)_n$ generated by the stochastic Douglas Rachford algorithm~\eqref{eq:sto-dr} with step $\gamma$, such that the distribution of $x_0^{\gamma,\nu}$ over $\sX$ is $\nu$. If $\nu = \delta_a$, where $\delta_a$ is the Dirac measure at the point $a \in \sX$, we shall prefer the notation $x_n^{\gamma,a}$. 
\section{Main convergence theorem}
\label{sec-th}

Consider the following assumptions.

\begin{assumption} 
\label{h0bnd} 
For every compact set $\cK\subset \sX$, there exists $\varepsilon>0$ such 
that
\[
\sup_{x\in\cK \cap \mD} 
\int \| \partial g_0(x,s) \|^{1+\varepsilon} \, \mu(ds) < \infty .   
\]
\end{assumption}

\begin{assumption}
\label{nablaf-bnd} 
For $\mu$-a.e $s \in \Xi$, $f(\cdot,s)$ is differentiable and there exists a closed ball in $\sX$ such that 
$\|\nabla f(x,s)\| \leq M(s)$ for all $x$ in this ball, where $M(s)$ is 
$\mu$-integrable. Moreover, for every compact set $\cK\subset \sX$, there 
exists $\varepsilon > 0$ such that 
\[
\sup_{x\in\cK} \int \| \nabla f (x,s) \|^{1+\varepsilon} \, \mu(ds) 
  < \infty\, . 
\] 
\end{assumption} 


\begin{assumption}
\label{linreg} 
$\displaystyle{
\forall x\in \sX, \ \int d(x,D(s))^2 \, \mu(ds) \geq C \bs d(x)^2}$.
\end{assumption}

\begin{assumption}
\label{JBnd-dif}
For every compact set $\cK\subset \sX$, there exists $\varepsilon, C, \gamma_0 > 0$ such that for all $\gamma \in (0,\gamma_0]$ and all $x \in \cK$, 
\[\frac 1{\gamma^{1+\varepsilon}}
\int \| \prox_{\gamma g(\cdot,s)}(x) - \Pi_{\cl(D(s))}(x) \|^{1+\varepsilon}  
 \, \mu(ds)  < C \, . 
\]
\end{assumption}

\begin{assumption}
\label{baillon}
There exists $L > 0$ such that $\nabla f(\cdot,s)$ is $\mu$-a.e, a $L$-Lipschitz continuous function.
\end{assumption}

\begin{assumption}
\label{L2} 
There exists $x_\star \in \arg\min F+G$ and $\varphi \in L^2(\Xi,\sX)$ such that $\varphi(s) \in \partial g(x_\star,s)$ $\mu$-a.s, $\nabla f(x_\star,\cdot) \in L^2(\Xi,\sX)$ and $\int \nabla f(x_\star,s) \, \mu(ds) + \int \varphi(s) \, \mu(ds) = 0$.
\end{assumption}

\begin{assumption}
\label{F+G-coerc} 
The function $F+G$ satisfies one of the following properties: 
\begin{enumerate}[label=(\alph*)]
\item\label{Zer-cpct} $F+G$ is coercive \textit{i.e} $F(x)+G(x) \longrightarrow_{\|x\| \to +\infty} +\infty$ 
\item\label{F+G-super} $F+G$ is supercoercive \textit{i.e} $\frac{F(x)+G(x)}{\|x\|} \longrightarrow_{\|x\| \to +\infty} +\infty$. 
\end{enumerate}
\end{assumption} 


\begin{assumption}
\label{JBgrow-fct}  
There exists $\gamma_0 > 0$, such that for all $\gamma\in (0,\gamma_0]$ and all $x\in \sX$, 
\begin{align*}
&\int \| \nabla f_\gamma(x,s) \| + \frac{1}{\gamma}\| \prox_{\gamma g(\cdot,s)}(x) - \Pi_{\cl(D(s))}(x) \|  \mu(ds)  
 \\ &\leq C ( 1 + | F^\gamma(x) + G^\gamma(x)| ) \, .  
\end{align*}
\end{assumption}


\begin{theorem}
\label{th-cv}
Let Assumptions~\ref{h0bnd}--~\ref{JBgrow-fct} hold true. Then, for each probability measure $\nu$ over $\sX$ having a finite second moment, for any $\varepsilon > 0$,
\[
\limsup_{n\to\infty} 
\frac 1{n+1}\sum_{k=0}^n 
 \bP\left( d(x_k^{\gamma,\nu}, \arg\min(F+G))>\varepsilon\right)
 \xrightarrow[\gamma\to 0]{}0\,. 
\] 
Moreover, if Assumption~\ref{F+G-coerc}--\ref{F+G-super} holds true, then 
\begin{gather*} 
\limsup_{n\to\infty}\ 
\bP\left(d\left(\bar x_n^{\gamma,\nu} , \arg\min(F+G) \right)\geq 
      \varepsilon\right)\xrightarrow[\gamma\to 0]{}0, \ \text{and} \\ 
\limsup_{n\to\infty}\ d\left(\bE(\bar x_n^{\gamma,\nu}), \arg\min(F+G) \right)
      \xrightarrow[\gamma\to 0]{}0\,. 
\end{gather*}
where $\bar x_n^{\gamma,\nu} = \frac{1}{n}\sum_{k=1}^{n} x_k^{\gamma,\nu}$.
\end{theorem} 

Loosely speaking, the theorem states that, with high probability, the iterates $(x_n^{\gamma,\nu})_n$ stay close to the set of solutions $\arg\min F+G$ as $n \to \infty$ and $\gamma \to 0$.  

Some Assumptions deserve comments. 

Following~\cite{bauschke1996projection}, we say that a finite collection of subsets $\mC_1,\dots,\mC_m$ of $\sX$ is \textit{linearly regular} if 
\begin{equation*}
\exists \kappa > 0, \forall x \in \sX, \max_{s \in \{1,\dots,m\}} d(x,\mC_s) \geq \kappa d(x,\cap_{s = 1}^{m} \mC_s)
\end{equation*}
In the case where there exists a $\mu$-probability one set $\tilde{\Xi}$ such that the set $\{D(s), s \in \tilde{\Xi}\} = \{\mC_1,\dots,\mC_m\}$ is finite, it is routine to check that Assumption~\ref{linreg} holds if and only if the domains $\mC_1,\dots,\mC_m$ are linearly regular. See~\cite{mourya2017adaptive} for an applicative context of the algorithm~\eqref{eq:sto-dr} in the latter case. 

It is a known fact that $$\prox_{\gamma g(\cdot,s)}(x) \longrightarrow_{\gamma \to 0} \Pi_{\cl({\dom(g(\cdot,s))})}(x),$$ for each $(x,s)$. Assumptions~\ref{JBnd-dif} and~\ref{JBgrow-fct} add controls on the convergence rate. 

Since $f(\cdot,s), g(\cdot,s) \in \Gamma_0(\sX)$, and $f(\cdot,s)$ is differentiable, $\partial (F+G)(x) = \nabla F(x) + \partial G(x) = \bE (\nabla f(x,\xi)) + \bE (\partial g(x,\xi))$ \cite{roc-wets-livre98}, where the set $\bE (\partial g(x,\xi))$ is defined by its Aumann integral $$\left\{\int \varphi(s) \, \mu(ds), \varphi \in L^1(\Xi,\sX), \text{ s.t. } \varphi(s) \in \partial g(x,s), \mu\text{-a.s.} \right\}$$
Therefore, using Fermat's rule, if $x \in \arg\min F+G$, then there exists $\varphi \in L^1(\Xi,\sX)$, such that $\varphi(s) \in \partial g(x,s)$ $\mu$-a.s, and $\int \nabla f(x,s) \, \mu(ds) + \int \varphi(s) \, \mu(ds) = 0$. We refer to $(\nabla f(x,\cdot), \varphi)$ as a \textit{representation} of the solution $x$. Assumption~\ref{L2} ensures the existence of $x_\star \in \arg\min F+G$ with a representation $\nabla f(x,\cdot), \varphi \in L^2(\Xi,\sX)$.



\section{Outline of the convergence proof}
\label{sec-proof}
This section is devoted to sketching the proof of the convergence of the stochastic Douglas Rachford algorithm. The approach follows the same steps as~\cite{bia-hac-sal-(sub)jca17} and is detailed in~\cite{salimSDR2017}.
The first step of the proof is to study the dynamical behavior of the iterates $(x_n^{\gamma,a})_n$ where $a \in \cD$. The Ordinary Differential Equation (ODE) method, well known in the literature of stochastic approximation (\cite{kus-yin-(livre)03}), is applied. Consider the continuous time stochastic process $\sx_{\gamma,a}$ obtained by linearly interpolating with time interval $\gamma$ the iterates $(x_n^{\gamma,a})$: 
\begin{equation}
\label{eq:interpol}
\sx_{\gamma,a}(t) = x_n^{\gamma,a} + (t - n\gamma)\frac{x_{n+1}^{\gamma,a} - x_n^{\gamma,a}}{\gamma},
\end{equation}
for all $t \geq 0$ such that $n\gamma \leq t < (n+1)\gamma$, for all $n \in \bN$. Let Assumptions~\ref{h0bnd}--\ref{JBnd-dif} \footnote{In the case where the domains are common, \textit{i.e} $s \mapsto D(s)$ is $\mu$-a.s constant, the moment Assumptions~\ref{h0bnd} and~\ref{nablaf-bnd} are sufficient to state the dynamical behavior result. See~\cite{mourya2017adaptive} for an applicative context where the domains $D(s)$ are distinct.} hold true. Consider the set $C(\bR_+,\sX)$ of continuous functions from $\bR_+$ to $\sX$ equipped
with the topology of uniform convergence on the compact intervals. It is shown that the continuous time stochastic process $\sx_{\gamma,a}$ converges weakly over $\bR_+$ (\textit{i.e} in distribution in $C(\bR_+,\sX)$) as $\gamma \to 0$. Moreover, the limit is proven to be the unique absolutely continuous function $\sx$ over $\bR_+$ satisfying $\sx(0) = a$ and for almost every $t \geq 0$, the Differential Inclusion (DI),  
\begin{equation}
\label{eq:di}
\dot \sx(t) \in - (\nabla F + \partial G)(\sx(t)), 
\end{equation}
(see \cite{bre-livre73}).
Differential inclusions like~\eqref{eq:di} generalize ODE to set-valued mappings. The DI~\eqref{eq:di} induces a map $\Phi : \mD \times \bR_+ \to \mD, (x_0, t) \mapsto \sx(t)$ that can be extended to a semi-flow over $\cl(\mD)$, still denoted by $\Phi$.

The weak convergence of $(\sx_{\gamma,a})$ to $\sx$ is not enough to study the long term behavior of the iterates $(x_n^{\gamma,a})_n$. The second step of the proof is to prove a stability result for the Feller Markov chain $(x_n^{\gamma,a})_n$. Denote by $P_\gamma$ its transition kernel. The deterministic counterpart of this step of the proof is the so-called \textit{Fejér monotonicity} of the sequence $(x_n)$ of the algorithm~\eqref{eq:dr-deter}. Even if some work has been done~\cite{bia-hac-jota16,combettes2015stochastic}, there is no immediate way to adapt the Fejér monotonicity to our random setting, mainly because of the constant step $\gamma$. As an alternative, we assume Hypotheses~\ref{baillon}-\ref{L2}, and prove the existence of positive numbers $\alpha, C$ and $\gamma_0$, such that for every $\gamma \in (0,\gamma_0]$,
\begin{align}
\label{eq:F+G}
\bE_{n} \|x_{n+1}^{\gamma,a}-x_\star\|^2 \leq &\| x_n^{\gamma,a} - x_\star\|^2\\
&- \alpha \gamma (F^\gamma + G^\gamma)(x_n^{\gamma,a}) + \gamma C. \nonumber 
\end{align}
In this inequality, $\bE_{n}$ denotes the conditional expectation with respect to the sigma-algebra $\sigma(x_0^\gamma,x_1^\gamma,\dots,x_n^\gamma)$ and $$F^\gamma(x) = \int f_\gamma(x,s) \, \mu(ds), \quad G^\gamma(x) = \int g_\gamma(x,s) \, \mu(ds).$$

Since $\gamma \mapsto F^\gamma(x) + G^\gamma(x)$ is decreasing~\cite{bia-hac-sal-(sub)jca17,salimSDR2017}, the function $F^\gamma+G^\gamma$ can be replaced by $F^{\gamma_0} + G^{\gamma_0}$. Besides, the coercivity of $F+G$ (Assumption~\ref{F+G-coerc}) implies the coercivity of $F^{\gamma_0} + G^{\gamma_0}$ (~\cite{bia-hac-sal-(sub)jca17,salimSDR2017}). Therefore, assuming~\ref{baillon}--\ref{F+G-coerc} and setting $\Psi = F^{\gamma_0} + G^{\gamma_0}$, there exist positive numbers $\alpha, C$ and $\gamma_0$, such that for every $\gamma \in (0,\gamma_0]$,
\begin{equation}
\label{eq:PH}
\bE_{n} \|x_{n+1}^{\gamma,a}-x_\star\|^2 \leq \| x_n^{\gamma,a} - x_\star\|^2
- \alpha \gamma \Psi (x_n^{\gamma,a}) + \gamma C. 
\end{equation}

Equation~\eqref{eq:PH} can alternatively be seen as a tightness result. 
It implies that the set $I_\gamma$ of invariant measures of the Markov kernel $P_\gamma$ is not empty for every $\gamma \in (0,\gamma_0]$, and that the set \begin{equation}\label{eq:Inv}
\text{Inv} = \cup_{\gamma \in (0,\gamma_0]} I_\gamma
\end{equation}
is \textit{tight}(~\cite{for-pag-99,bia-hac-sal-(arxiv)16}).

It remains to characterize the cluster points of Inv as $\gamma \to 0$. To that end, the dynamical behavior result and the stability result are combined. Let Assumptions~\ref{h0bnd}--~\ref{JBgrow-fct} hold true. \footnote{Assumptions~\ref{linreg},~\ref{JBnd-dif} and~\ref{JBgrow-fct} are not needed if the domains $D(s)$ are common.} Then, the set Inv is tight, and, as $\gamma \to 0$, every cluster point of Inv is an invariant measure for the semi-flow $\Phi$. The Theorem~\ref{th-cv} is a consequence of this fact.



\section{Application to structured regularization}
\label{sec-simu}
In this section is provided an application of the stochastic Douglas Rachford~\eqref{eq:sto-dr} algorithm to solve a regularized optimization problem. Consider problem~\eqref{eq:minF+G}, where $F$ is a cost function that is written as an expectation, and $G$ is a regularization term. Towards solving~\eqref{eq:minF+G}, many approaches involve the computation of the proximity operator of the regularization term $G$. In the case where $G$ is a structured regularization term, its proximity operator is often difficult to compute. When $G$ is a graph-based regularization, it is possible to apply a stochastic proximal method to address the regularization~\cite{sal-bia-hac-(sub)tac17}. We shall concentrate on the case where $G$ is an overlapping group regularization. In this case, the computation of the proximity operator of $G$ is known to be a bottleneck~\cite{yuan2011efficient}. We shall apply the algorithm~\eqref{eq:sto-dr} to overcome this difficulty.

Consider $\sX = \bR^N$, $N \in \bN^\star$, and $g \in \bN^\star$. Consider $g$ subsets of $\{1,\dots,N\}$, $S_1,\dots,S_g$, possibly overlapping. Set $G(x) = \sum_{j = 1}^{g} \|x_{S_j}\|$, where $x_{S_j}$ denotes the restriction of $x$ to the set of index $S_j$ and $\|\cdot\|$ denotes the Euclidean norm. Set $F(x) = \bE_{(\xi,\eta)}(h(\eta \ps{x,\xi}))$ where $h$ denotes the hinge loss $h(z) = \max(0,1-z)$ and $(\xi,\eta)$ is a r.v defined on some probability space with values in $\sX \times \{-1,+1\}$. In this case, the problem~\eqref{eq:minF+G} is also called the SVM classification problem, regularized by the overlapping group lasso. It is assumed that the user is provided with i.i.d copies $((\xi_n,\eta_n))_n$ of the r.v $(\xi,\eta)$ online. 

To solve this problem, we implement a stochastic Douglas Rachford strategy. To that end, the regularization $G$ is rewritten $G(x) = \bE_J(g \|x_{S_J}\|)$ where $J$ is an uniform r.v over $\{1,\dots,g\}$. At each iteration $n$ of the stochastic Douglas Rachford algorithm, the user is provided with the realization $(\xi_n,\eta_n)$ and sample a group $J_n$ uniformly in $\{1,\dots,g\}$. Then, a Douglas Rachford step is done, involving the computation of the proximity operators of the functions $g_n : x \mapsto \|x_{S_{J_n}}\|$ and $f_n : x \mapsto h(\eta_n \ps{x,\xi_n})$. 

\begin{figure}[ht!]
 \includegraphics[width=\linewidth]{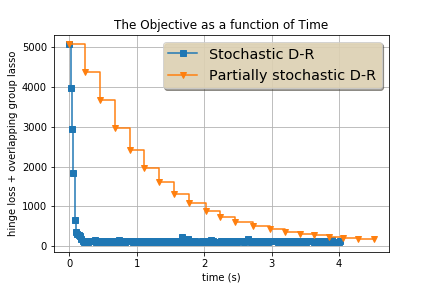}
\caption{The objective function $F+G$ as a function
of time in seconds for each algorithm}
\label{fig:simu}
\end{figure}

This strategy is compared with a partially stochastic Douglas Rachford algorithm, deterministic in the regularization $G$, where the fast subroutine Fog-Lasso~\cite{yuan2011efficient} is used to compute the proximity operator of the regularization $G$. At each iteration $n$, the user is provided with $(\xi_n,\eta_n)$. Then, a Douglas Rachford step is done, involving the computation of the proximity operators of the functions $G$ and $f_n : x \mapsto h(\eta_n \ps{x,\xi_n})$. Figure~\ref{fig:simu} demonstrates the advantage of treating the regularization term in a stochastic way. 

\begin{figure}[ht!]
 \includegraphics[width=\linewidth]{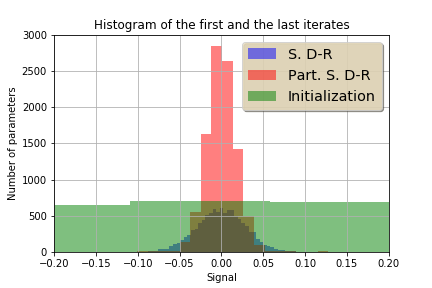}
\caption{Histogram of the Initialization and the last iterates of the Stochastic D-R (S. D-R) and the partially stochastic D-R (Part. S. D-R)}
\label{fig:hist}
\end{figure}

In Figure~\ref{fig:simu} "Stochastic D-R" denotes the stochastic Douglas Rachford algorithm and "Partially stochastic D-R" denotes the partially stochastic Douglas Rachford where the subroutine FoG-Lasso~\cite{yuan2011efficient} is used at each iteration to compute the true proximity operator of the regularization $G$. Figure~\ref{fig:hist} shows the appearance of the first and the last iterates. Even if a best performing procedure~\cite{yuan2011efficient} is used to compute $\prox_{\gamma G}$, we observe on Figure~\ref{fig:simu} that Stochastic D-R takes advantage of being a stochastic method. This advantage is known to be twofold (\cite{bottou2016optimization}). First, the iteration complexity of Stochastic D-R is
moderate because $\prox_{\gamma G}$ is never computed. Then, Stochastic D-R is faster than its partially deterministic counterpart which uses Fog-Lasso~\cite{yuan2011efficient} as a subroutine, especially in the first
iterations of the algorithms.
Moreover, Stochastic D-R seems to perform globally better. This
is because every proximity
operators in Stochastic D-R can be efficiently computed (\cite{bau-com-livre11}). Contrary to the proximity operator of $G$~\cite{yuan2011efficient}, the proximity operator of $g_n$ is easily computable. The proximity operator of $f_n$ is easily computable as well.\footnote{Even if $h(x) = \log(1+\exp(-x))$ (logistic regression), the proximity operator of $f_n$ is easily computable, see~\cite{chierchiaapproche}.}

\bibliographystyle{IEEEbib}
\bibliography{math}

\end{document}